\documentclass[a4paper,10pt]{article}
\usepackage[T1]{fontenc}
\usepackage[utf8]{inputenc}
\usepackage{amsmath,amssymb}
\usepackage{mathtools} 
\usepackage{subcaption}
\usepackage{graphicx}
\usepackage{enumitem}
\setenumerate{label=\alph*)}
\usepackage{bm}
\usepackage{authblk}
\usepackage{amsthm}
\usepackage[margin=1.2in]{geometry}

\usepackage{graphicx}

\newcommand{\bn}{\mathbf n}

\newcommand{\bg}{\mathbf g}
\newcommand{\bA}{\mathbf A}

\newcommand{\bN}{\mathbf N}
\newcommand{\bB}{\mathbf B}

\newcommand{\bD}{\mathbf D}

\newcommand{\bI}{\mathbf I}

\newcommand{\bff}{\mathbf f}

\newcommand{\bzero}{\mathbf 0}
\newcommand{\bv}{\mathbf v}

\newcommand{\by}{\mathbf y}

\newcommand{\R}{\mathbb R}

\newcommand{\tm}{\subseteq}
\newcommand{\abs}[1]{\lvert #1\rvert}
\newcommand{\norm}[1]{\lVert #1\rVert}

\renewcommand{\R}{\mathbb{R}}
\newcommand{\N}{\mathbb{N}}

\newcommand{\bby}{\bar{\by}}
\renewcommand{\Vec}[1]{\renewcommand*{\arraystretch}{1}\begin{pmatrix*}[r]#1\end{pmatrix*}}
\renewcommand{\vec}[1]{\begin{psmallmatrix*}[r]#1\end{psmallmatrix*}}

\newcommand{\from}{\colon}

\DeclareMathOperator{\Span}{span}

\makeatletter
\newcommand{\vast}{\bBigg@{3}}
\newcommand{\Vast}{\bBigg@{4}}
\makeatother

\theoremstyle{definition}
\newtheorem{thm}{Theorem}[section]
\newtheorem{defn}[thm]{Definition}

\newtheorem{rem}[thm]{Remark}
\newtheorem{cor}[thm]{Corollary}

\begin{document}

		\title{A Stability Analysis of Modified
			Patankar--Runge--Kutta methods for a nonlinear Production--Destruction System }
		
\author[1]{Thomas Izgin} 
\author[1]{Stefan Kopecz}
\author[1]{Andreas Meister}
\affil[1]{Department of Mathematics, University of Kassel, Germany}
\affil[]{izgin@mathematik.uni-kassel.de\ \&\ kopecz@mathematik.uni-kassel.de\ \&\ meister@mathematik.uni-kassel.de}
\setcounter{Maxaffil}{0}
\renewcommand\Affilfont{\itshape\small}
\date{}
		\maketitle  
		\begin{abstract}
		Modified Patankar--Runge--Kutta (MPRK) methods preserve the positivity as well as conservativity of a production--destruction system (PDS) of ordinary differential equations for all time step sizes. As a result, higher order MPRK schemes do not belong to the class of general linear methods, i.\,e.\ the iterates are generated by a nonlinear map $\mathbf g$ even when the PDS is linear. Moreover, due to the conservativity of the method, the map $\mathbf g$ possesses non-hyperbolic fixed points.
	
		Recently, a new theorem for the investigation of stability properties of non-hyperbolic fixed points of a nonlinear iteration map was developed. We apply this theorem to understand the stability properties of a family of second order MPRK methods when applied to a nonlinear PDS of ordinary differential equations. It is shown that the fixed points are stable for all time step sizes and members of the MPRK family. Finally, experiments are presented to numerically support the theoretical claims.
		\end{abstract}
		
		\section{Introduction}
The mathematical modeling of numerous applications of natural sciences and engineering leads to positive and conservative production-destruction systems (PDS) 
\[
y_i'(t)=\sum_{j=1}^N( p_{ij}(\mathbf y(t))-d_{ij}(\mathbf y(t))),\quad\by(0)=\by^0\in\R^N_{>0},\quad p_{ij}(\mathbf y),d_{ij}(\mathbf y)\geq 0,\quad i=1,\dotsc, N,
\]
where $\mathbf y=(y_1,\dots,y_N)^T$ denotes the vector of state variables. The terms $p_{ij}, d_{ij}$ denote production and destruction terms of the $i$-th constituent, respectively. A PDS is called \emph{conservative}, if $p_{ij}=d_{ji}$, i.\,e.\ $\sum_{i=1}^N y_i(t)=\sum_{i=1}^N y_i(0)$ is satisfied for all $t\geq 0$. The PDS is called \emph{positive}, if $\by(t)>\bzero$ holds for all $t> 0$ whenever $\by(0)>\bzero$. 

When applied to a positive and conservative PDS, a numerical method should satisfy these properties in a discrete manner. This means that we call the method  \emph{unconditionally positive} if $\by^n$ is positive for all $\by^0>\bzero$, $n\in \N$ and $\Delta t>0$, and \emph{unconditionally conservative} if $\sum_{i=1}^N y_i^{n+1}=\sum_{i=1}^N y_i^{n}$ holds true for all $n\in \N_0$ and $\Delta t>0$.

One class of unconditionally positive and conservative numerical methods of second order is given by the two-stage one-parameter family of modified Patankar--Runge--Kutta methods denoted by MPRK22($\alpha$). This family of schemes is based on the set of explicit two-stage Runge--Kutta (RK) methods with nonnegative parameters given by the Butcher tableau
$$
\begin{array}{c|cc}
	0        &                    & \\
	\alpha & \alpha            & \\\hline
	& 1-1/(2\alpha) & 1/(2\alpha)
\end{array},\quad \alpha\geq \frac12, 
$$
and defined by
\begin{subequations}\label{eq:MPRK22b}
	\begin{align*}
		&\begin{aligned}
			\mathllap{y_i^{(1)}} &= y_i^n,\quad\quad\quad \mathllap{y_i^{(2)}} = y_i^n + \alpha\Delta t\sum_{j=1}^N\biggl( p_{ij}(\by^{(1)})\frac{y_j^{(2)}}{y_j^{(1)}}-d_{ij}(\by^{(1)})\frac{y_i^{(2)}}{y_i^{(1)}}\biggr),
		\end{aligned}\\ 
		&\begin{multlined}[b][.7\columnwidth]
			\mathllap{y_i^{n+1}} = y_i^n + \Delta t\sum_{j=1}^N\biggl( \left(\left(1-\frac1{2\alpha}\right) p_{ij}(\by^{(1)})+\frac1{2\alpha} p_{ij}(\by^{(2)})\right)\frac{y_j^{n+1}}{(y_j^{(2)})^{1/\alpha}(y_j^{(1)})^{1-1/\alpha}}\\
			- \left(\left(1-\frac1{2\alpha}\right) d_{ij}(\by^{(1)})+ \frac1{2\alpha} d_{ij}(\by^{(2)})\right)\frac{y_i^{n+1}}{(y_i^{(2)})^{1/\alpha}(y_i^{(1)})^{1-1/\alpha}}\biggr),
		\end{multlined}  
	\end{align*}
\end{subequations}
for $i=1,\dots,N$, see \cite{KM18}. 

Since the numerical method should also be capable of replicating the stability properties of the underlying problem, 
steady states of the differential equation should be fixed points of the numerical scheme with identical stability properties.

\begin{defn}\label{Def Lyapunov Cont}
	Let $\by^*\in \R^N$ be a steady state solution of a differential equation $\by'={\bff}(\by)$, that is ${\bff}(\by^*)=\bzero$.
	\begin{enumerate}
		\item\label{item:lyap_stab} Then $\by^*$ is called \emph{Lyapunov stable} if, for any $\epsilon>0$, there exists a $\delta=\delta(\epsilon)>0$ such that $\norm{\by(0)-\by^*}<\delta$ implies $\norm{\by(t)-\by^*}<\epsilon$ for all $t\geq 0$.
		\item If in addition to a), there exists a constant $c>0$ such that $\Vert \by(0)-\by^*\Vert<c$ implies  $\Vert \by(t)-\by^*\Vert \to 0$ for $t\to \infty$,  we call  $\by^*$ \emph{asymptotically stable.}
		\item A steady state solution that is not Lyapunov stable is said to be \emph{unstable}.
	\end{enumerate}
\end{defn}
In the following we will also briefly speak of stability instead of Lyapunov stability.

With these notions in mind let us consider the nonlinear test equation	
\begin{equation}
	\by'(t)=\bff(\by)=\vec{
		y_2^2-y_1^2 \\ y_1^2-y_2^2}\label{eq:PDS_test}
\end{equation}
together with the initial condition 
\begin{equation}
	\by(0)=\by^0 = \vec{
		y_1^0 \\ y_2^0} > \bzero.\label{eq:IC}
\end{equation}
The test problem \eqref{eq:PDS_test} can be written as a PDS using
\begin{equation}\label{eq:pij,dij}
	p_{ij}(\by)=d_{ji}(\by)=y_j^2
\end{equation}
for $i,j=1,2$ with $i\neq j$ and $p_{ii}(\by)=d_{ii}(\by)=0$ for $i=1,2$.
Using $\bm 1=(1,1)^\intercal$, the set of positive steady states of \eqref{eq:PDS_test} is given by $\Span(\bm 1)\cap \R^2_{>0}$, and the solution of the associated initial value problem given by \eqref{eq:PDS_test} and \eqref{eq:IC} can be written as
\begin{equation}\label{eq:anasol}
	\by(t)=\frac12 (y_1^0+y_2^0)\bm 1+\frac12\Bigl(\by^0-\vec{y_2^0\\y_1^0}\Bigr) e^{-2(y_1^0 + y_2^0)t}.
\end{equation}
Due to $\by^0>\bzero$, the exponential term vanishes as $t\to \infty$ and hence, $\lim_{t\to\infty}\by(t)=\frac12 (y_1^0+y_2^0)\bm 1\in \Span(\bm 1)\cap \R^2_{>0}$. In particular, given a positive initial condition, the exact solution monotonically approaches the positive steady state solution
\begin{equation}
	\by^*=\frac12 (y_1^0+y_2^0)\bm 1\label{eq:y*}
\end{equation}
along the line $y_1+y_2=y_1^0+y_2^0$ in the $y_1$-$y_2$-coordinate system. As a result, $\norm{\by(t)-\by^*}<\norm{\by(0)-\by^*}$ is true for all $t>0$, which means that we can choose $\delta=\epsilon$ in Definition \ref{Def Lyapunov Cont} to see that all positive steady states are stable. However, none of them is asymptotically stable since in any neighborhood of a steady state there are infinitely many further steady states.

In total the numerical method should transfer stable but not asymptotically stable steady states to fixed points with similar properties, which we define analogously to the continuous case. 
\begin{defn}\label{Def_Lyapunov_Diskr}
	Let $\by^*\in \R^N$ be a fixed point of an iteration scheme $\by^{n+1}=\bg(\by^n)$, that is $\by^*=\bg(\by^*)$. 
	\begin{enumerate}
		\item\label{def:stab} Then $\by^*$ is called \emph{Lyapunov stable} if, for any $\epsilon>0$, there exists a $\delta=\delta(\epsilon)>0$ such that $\norm{\by^0-\by^*}<\delta$ implies $\norm{\by^n- \by^*}<\epsilon$ for all $n\geq 0$.
		\item If in addition to a), there exists a constant $c>0$ such that $\Vert \by^0-\by^*\Vert<c$ implies $\Vert \by^n-\by^*\Vert \to 0$ for $n\to \infty$, we call $\by^*$ \emph{asymptotically stable.}
		\item A fixed point that is not Lyapunov stable is said to be \emph{unstable}.
	\end{enumerate}
\end{defn}
Analogous to steady states, we also speak of stability instead of Lyapunov stability in the case of fixed points.

So far, the stability of MPRK schemes, see \cite{KM18,KM18b,HS,HZS,OeT} has been studied exclusively in the context of linear systems of differential equations in \cite{IKM22,IKM2122,HIKMS22,IOE22StabMP}. Here we want to show that the same approach can be used to investigate the stability of MPRK schemes applied to the nonlinear problem \eqref{eq:PDS_test}.

To recap the stability theorem from \cite{IKM22}, we introduce a matrix $\bA\in \R^{N\times N}$ such that $\bn_1,\dotsc,\bn_k$ with $k\geq 1$ form a basis of $\ker(\bA^\intercal)$.  We also define the matrix
\begin{equation*}\label{eq:N}
	\bN=\begin{pmatrix}
		\bn_1^\intercal\\
		\vdots\\
		\bn_k^\intercal
	\end{pmatrix}\in \R^{k\times N} 
\end{equation*}
and the set
\begin{equation}
	H=\{\by\in \R^N\mid \bN\by=\bN\by^*\}\label{eq:H}.
\end{equation}
With this in mind, we can now state the following theorem and point out that it is a generalization of \cite[Theorem 2.9]{IKM2122} to systems of arbitrary finite size. 
\begin{thm}[{\cite[Theorem 2.9]{IKM22}}]\label{Thm_MPRK_stabil}
	Let $\bA\in \R^{N\times N}$ such that $\ker(\bA)=\Span(\bv_1,\dotsc,\bv_k)$ represents a $k$-dimensional subspace of $\R^N$ with $k\geq 1$. Also, let $\by^*\in \ker(\bA)$ be a fixed point of $\bg\from D\to D$ where $D\tm \R^N$ contains a neighborhood $\mathcal D$ of $\by^*$. Moreover, let any element of $\ker(\bA)\cap \mathcal D$ be a fixed point of $\bg$ and suppose that $\bg\big|_\mathcal{D}\in \mathcal C^1$ as well as that the first derivatives of $\bg$ are Lipschitz continuous on $\mathcal{D}$. Then $\bD\bg(\by^*)\bv_i=\bv_i$ for $i=1,\dotsc, k$ and the following statements hold.
	\begin{enumerate}
		\item\label{it:Thma} If the remaining $N-k$ eigenvalues of $\bD\bg(\by^*)$ have absolute values smaller than $1$, then $\by^*$ is stable.\label{It:Thm_Stab_a}
		\item\label{it:Thmb} Let $H$ be defined by \eqref{eq:H} and $\bg$ conserve all linear invariants, which means that $\bg(\by)\in H\cap D$ for all $\by\in H\cap D$. If additionally the assumption of \ref{It:Thm_Stab_a} is satisfied, then there exists a $\delta>0$ such that $\by^0\in H\cap D$ and $\norm{\by^0-\by^*}<\delta$ imply $\by^n\to \by^*$ as $n\to \infty$.
	\end{enumerate}
\end{thm}
\begin{rem}
	In order to use Theorem \ref{Thm_MPRK_stabil} to analyze MPRK22($\alpha$) applied to \eqref{eq:PDS_test} we note that
	\begin{enumerate}
		\item  the matrix $\bA=\Vec{1 &-1\\-1 & 1}$ satisfies $\ker(\bA)=\ker(\bA^\intercal)=\Span(\bm 1)\tm \{\by\in \R^N\mid \bff(\by)=\bzero\}$,
		\item any positive steady state of \eqref{eq:PDS_test} is a fixed point of MPRK22($\alpha$), see \cite{issuesTOER}, 
		\item $\bg$ conserves all linear invariants as the only linear invariant is conservativity, and
		\item it is sufficient to prove $\bg\in \mathcal C^2$ rather than that $\bg\in\mathcal C^1$ has locally Lipschitz first derivatives, see \cite{IKM2122, IKM22}.
	\end{enumerate}
	As a result of \cite[Remark 2.10]{IKM22} it remains to prove that $\bg\in \mathcal C^2$ in a small enough neighborhood of a positive steady state $\by^*$ and to analyze the eigenvalues of the Jacobian $\bD\bg(\by^*)$.
\end{rem}

We first apply the MPRK22($\alpha$) method to \eqref{eq:PDS_test} 
obtaining
\begin{align*}
	y_i^{(2)}&=y_i^n+\alpha\Delta t (y_j^ny_j^{(2)}-y_i^ny_i^{(2)}),\\
	y_i^{n+1}&=y_i^n+\Delta t \left(\left(1-\tfrac{1}{2\alpha}\right)(y_j^n)^2+\frac{1}{2\alpha}(y_j^{(2)})^2\right)\frac{y_j^{n+1}}{(y_j^{(2)})^{\tfrac{1}{\alpha}}(y_j^{n})^{1-\tfrac{1}{\alpha}}}\\
	&\hphantom{=y_i^n}-\Delta t \left(\left(1-\tfrac{1}{2\alpha}\right)(y_i^n)^2+\frac{1}{2\alpha}(y_i^{(2)})^2\right)\frac{y_i^{n+1}}{(y_i^{(2)})^{\tfrac{1}{\alpha}}(y_i^{n})^{1-\tfrac{1}{\alpha}}}.
\end{align*}
Following the approach described in \cite{IKM22,HIKMS22,IOE22StabMP}, we define the functions $\bm \Phi_2$ and $\bm\Phi_{n+1}$ by
\begin{align}
	\bm \Phi_2(\by^n,\by^{(2)})_i&=-y_i^{(2)}+y_i^n+\alpha\Delta t (y_j^ny_j^{(2)}-y_i^ny_i^{(2)}),\label{eq:Phi2}\\
	\bm \Phi_{n+1}(\by^n,\by^{(2)},\by^{n+1})_i&=-	y_i^{n+1}+y_i^n+\Delta t \left(\left(1-\tfrac{1}{2\alpha}\right)(y_j^n)^2+\frac{1}{2\alpha}(y_j^{(2)})^2\right)\frac{y_j^{n+1}}{(y_j^{(2)})^{\tfrac{1}{\alpha}}(y_j^{n})^{1-\tfrac{1}{\alpha}}}\nonumber\\
	&\hphantom{=y_i^n}-\Delta t \left(\left(1-\tfrac{1}{2\alpha}\right)(y_i^n)^2+\frac{1}{2\alpha}(y_i^{(2)})^2\right)\frac{y_i^{n+1}}{(y_i^{(2)})^{\tfrac{1}{\alpha}}(y_i^{n})^{1-\tfrac{1}{\alpha}}}\label{eq:Phin+1}
\end{align}
for $i,j=1,2$ and $i\neq j$. The functions $\bm \Phi_2=\bm \Phi_2(\bm u, \bm v)$ and $\bm \Phi_{n+1}=\bm \Phi_{n+1}(\bm u, \bm v,\bm w)$ are in $\mathcal C^2$ for $\bm u,\bm v,\bm w\in\R^2_{>0}$, so that we can define 
\begin{equation*}
	\bD_{n}\bm\Phi_k=\frac{\partial\bm\Phi_k}{\partial\bm u},\quad \bD_{2}\bm\Phi_k=\frac{\partial\bm\Phi_k}{\partial\bm v} \quad \text{for  }k=2,n+1
\end{equation*}
as well as
\begin{equation*}
	\bD_{n+1}\bm\Phi_{n+1}=\frac{\partial\bm\Phi_{n+1}}{\partial\bm w}.
\end{equation*}
It is worth mentioning that the operator $\bD^*$ applied to $\bm\Phi_k$ for $k=2,n+1$ means that we plug in $\by^*$ for all arguments, e.\,g.\ \[\bD^*_n\bm\Phi_2=\bD_n\bm\Phi_2(\by^*,\by^*).\] If the inverses of $\bD^*_{n+1}\bm\Phi_{n+1}$ and $\bD_2\bm\Phi_2$ exist we can conclude 
\begin{equation}\label{eq:Dg(y*)}
	\bD\bg(\by^*)=-(\bD^*_{n+1}\bm\Phi_{n+1})^{-1}(\bD^*_n\bm \Phi_{n+1}-\bD^*_2\bm \Phi_{n+1}(\bD^*_2\bm \Phi_{2})^{-1}\bD^*_n\bm\Phi_2),
\end{equation}
see \cite{IOE22StabMP}. Moreover, in this case the implicit function theorem states that also $\bg\in \mathcal C^2$ in a small enough neighborhood of $\by^*$. As a result, all we have to do is prove that the inverses occurring in \eqref{eq:Dg(y*)} exist and to investigate the spectrum of $\bD\bg(\by^*)$.

Let us introduce the matrix \[\bB=\Vec{-1 & 1\\1 &-1},\]
so that \eqref{eq:Phi2} and $\by^*=c\bm 1$ with $c=\frac12(y_1^0+y_2^0)$ yield
\begin{equation*}
	\begin{aligned}
		\bD_n^*\bm\Phi_2=\Vec{1-\alpha\Delta tc & \alpha\Delta tc\\ \alpha\Delta tc & 1-\alpha\Delta tc }=\bI-\alpha\Delta tc\bB	
	\end{aligned}
\end{equation*}
and 
\begin{equation*}
	\bD_2^*\bm\Phi_2=\Vec{-1-\alpha\Delta tc & \alpha\Delta tc\\ \alpha\Delta tc & -1-\alpha\Delta tc }=-\bI+\alpha\Delta tc\bB.
\end{equation*}
Note that $\sigma(\bB)=\{0,-2\}$ and $\alpha\Delta tc>0$ imply that the inverse of $\bD_2^*\bm\Phi_2$ exists. Next, using \eqref{eq:Phin+1} and the quotient rule, we obtain the Jacobian
\begin{equation}\label{eq:DnPhin+1}
	\begin{aligned}
		\bD_n^*\bm\Phi_{n+1}&=\Vec{1-\Delta t \frac{2(1-\tfrac{1}{2\alpha})(y_1^*)^3-(y_1^*)^3(1-\tfrac{1}{\alpha})}{(y_1^*)^2} & \Delta t(2(1-\tfrac{1}{2\alpha})y_2^*-y_2^*(1-\tfrac{1}{\alpha})\\\Delta t(2(1-\tfrac{1}{2\alpha})y_1^*-y_1^*(1-\tfrac{1}{\alpha}) & 1-\Delta t \frac{2(1-\tfrac{1}{2\alpha})(y_2^*)^3-(y_2^*)^3(1-\tfrac{1}{\alpha})}{(y_2^*)^2}}\\
		&=\Vec{1-c\Delta t & c\Delta t\\ c\Delta t & 1-c\Delta t}=\bI+c\Delta t\bB.
	\end{aligned}
\end{equation}
Furthermore, we have
\begin{equation}\label{eq:D2Pn+1}
	\begin{aligned}
		\bD_2^*\bm\Phi_{n+1}&=\Vec{\Delta t \frac{\tfrac{1}{\alpha}(y_1^*)^3-(y_1^*)^3\tfrac{1}{\alpha}}{(y_1^*)^2} & \Delta t(\tfrac{1}{\alpha}y_2^*-y_2^*\tfrac{1}{\alpha})\\\Delta t(\tfrac{1}{\alpha}y_1^*-y_1^*\tfrac{1}{\alpha}) & \frac{\tfrac{1}{\alpha}(y_2^*)^3-(y_2^*)^3\tfrac{1}{\alpha}}{(y_2^*)^2}}=\bzero
	\end{aligned}
\end{equation}
as well as the matrix
\begin{equation}\label{eq:Dn+1Phin+1}
	\begin{aligned}
		\bD_{n+1}^*\bm\Phi_{n+1}&=\Vec{-1-\Delta ty_1^*  & \Delta t y_2^*\\ \Delta ty_1^* & -1-\Delta ty_2^*}=-\bI+\Delta tc\bB,
	\end{aligned}
\end{equation}
which is nonsingular due to $\Delta tc \geq 0$. Altogether, we are able to use the formula \eqref{eq:Dg(y*)}, and due to \eqref{eq:D2Pn+1} it simplifies to 
\begin{equation*}
	\bD\bg(\by^*)=-(\bD^*_{n+1}\bm\Phi_{n+1})^{-1}\bD^*_n\bm \Phi_{n+1}.
\end{equation*}
Using \eqref{eq:DnPhin+1} and \eqref{eq:Dn+1Phin+1} we obtain
\begin{equation*}
	\begin{aligned}
		\bD\bg(\by^*)&=(\bI-\Delta tc\bB)^{-1}(\bI+\Delta tc\bB)=\frac{1}{2\Delta tc+1}\Vec{1+\Delta tc &\Delta tc\\ \Delta tc & 1+\Delta tc}\Vec{1- \Delta tc & \Delta tc\\\Delta tc &1-\Delta tc}\\
		&=\frac{1}{2\Delta tc+1}\Vec{1 & 2\Delta tc\\ 2\Delta tc &1}.
	\end{aligned}
\end{equation*}
In accordance with \cite{IKM2122} we find that the eigenvectors are given by $\by^*$ and $\bby=(1,-1)^\intercal$ satisfying $\bD\bg(\by^*)\by^*=\by^*$ and 
\begin{equation*}
	\bD\bg(\by^*)\bby=\frac{1-2\Delta t c}{1+2\Delta tc}\bby.
\end{equation*} 
Hence, using Theorem \ref{Thm_MPRK_stabil} and substituting $z=2\Delta tc$, we have to analyze the stability function 
\begin{equation}
	R(z)=\frac{1-z}{1+z}.
\end{equation}
Indeed, $R(0)=1$, $\lim_{z\to\infty}R(z)=-1$ and $R'(z)=-\tfrac{2}{(1+z)^2}<0$ yield
\[\lvert R(z)\rvert <1 \quad \text{for all $z=2\Delta tc>0$,}\]
and hence, Theorem \ref{Thm_MPRK_stabil} implies the following results.
\begin{cor}\label{Cor:MPRKstab}
	Let $\by^*$ be a positive steady state of the differential equation \eqref{eq:PDS_test}. Then $\by^*$ is a stable fixed point of the MPRK22($\alpha$) scheme for all $\Delta t>0$ and $\alpha\geq \tfrac12$.
\end{cor}
\begin{cor}\label{Cor:MPRKstab1}
	Let the unique steady state $\by^*$ of the initial value problem \eqref{eq:PDS_test}, \eqref{eq:IC} be positive. Then the iterates of MPRK22($\alpha$) locally converge towards $\by^*$ for all $\Delta t>0$ and $\alpha\geq \tfrac12$.
\end{cor}
\section{Numerical Experiments}
\begin{figure}[!h]
	\centering
	\includegraphics[width=0.4\textwidth]{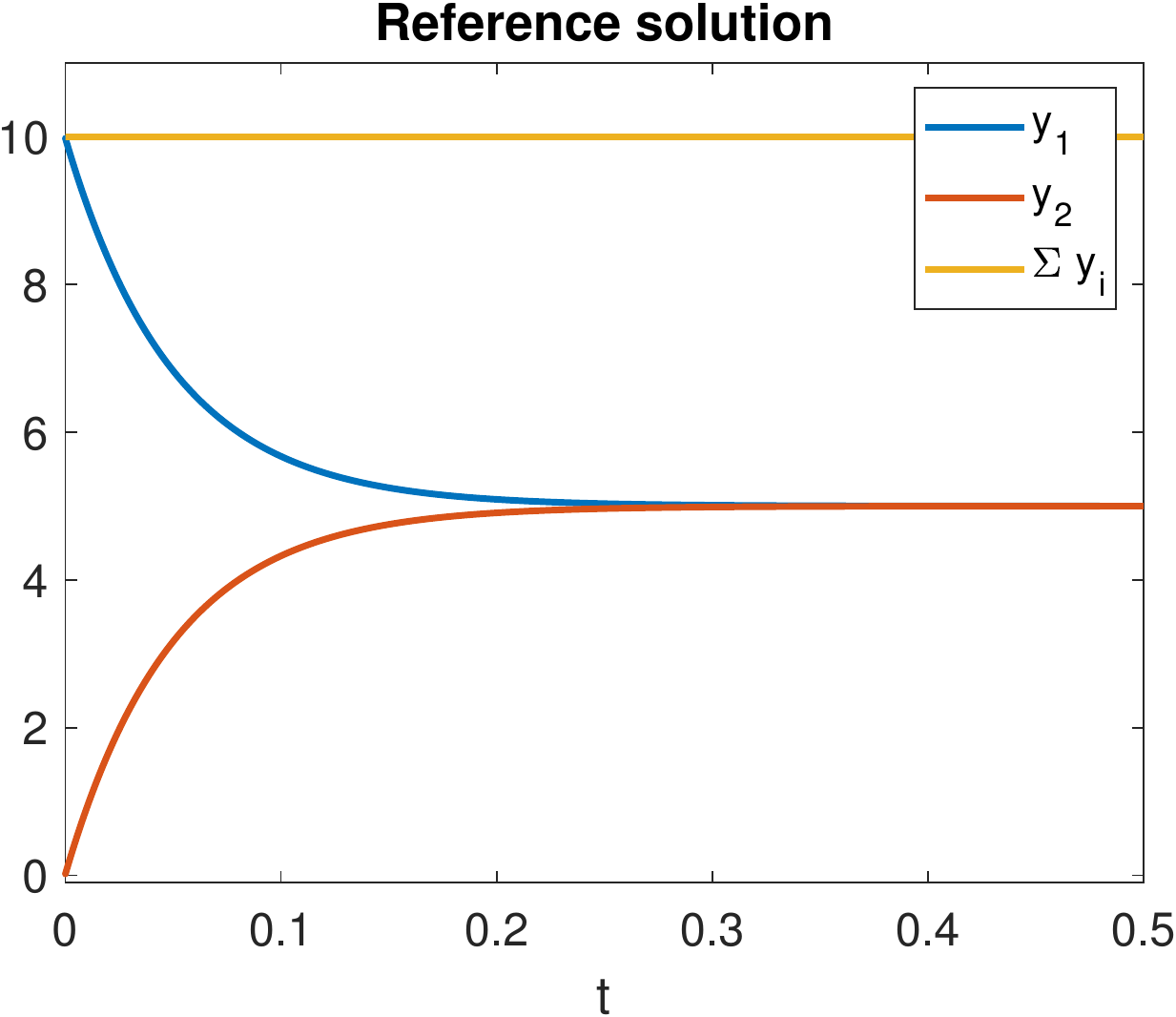}
	\caption{Analytical solution \eqref{eq:anasol} of \eqref{eq:PDS_test} with $\by^0=(9.98,0.02)^\intercal$.}\label{fig:analytic}
\end{figure}

\begin{figure}[!h]
	\centering
	\begin{subfigure}[t]{0.4\textwidth}
		\includegraphics[width=\textwidth]{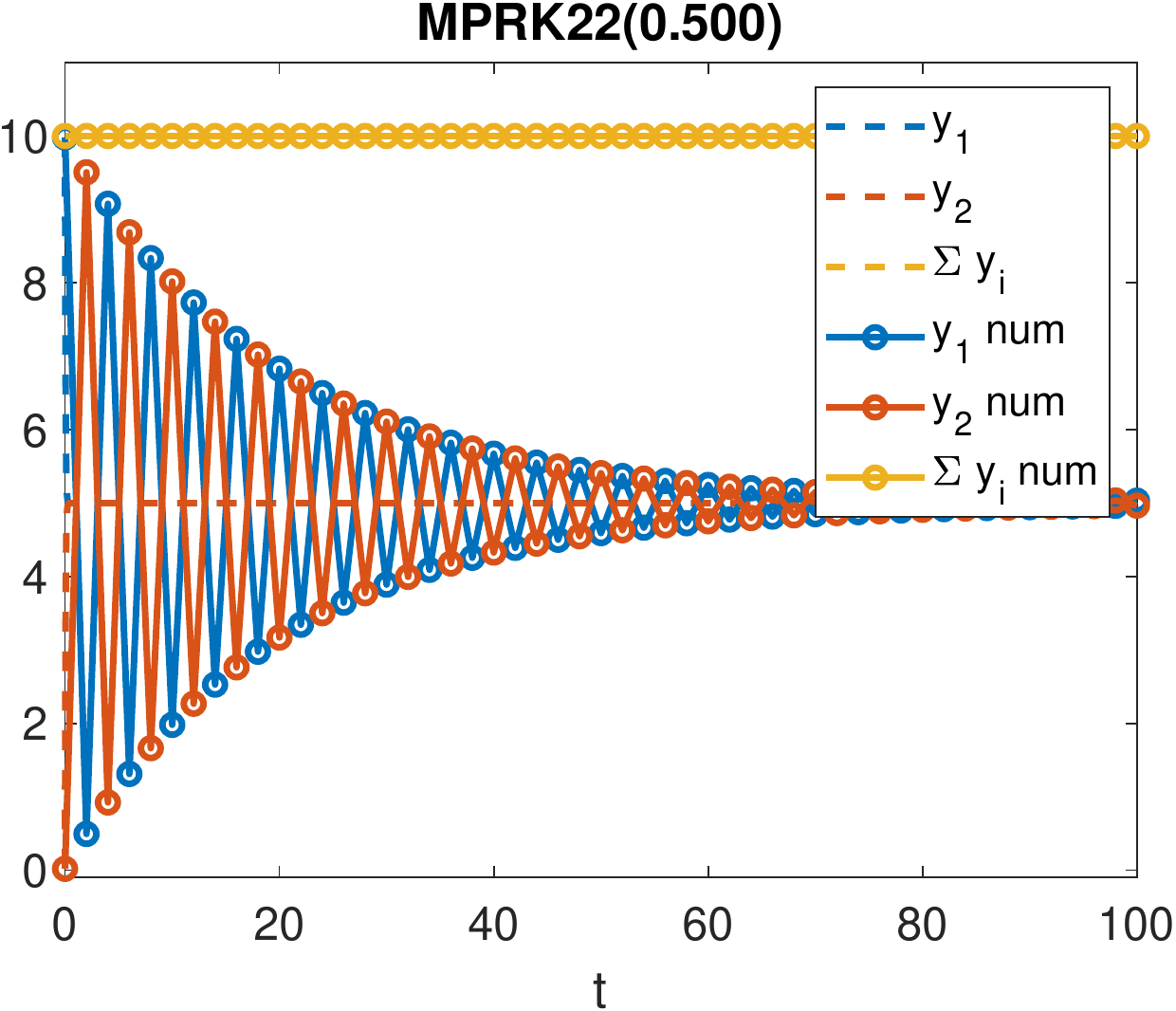}
		\subcaption{}\label{subfig:0.5}
	\end{subfigure}
	\begin{subfigure}[t]{0.423\textwidth}
		\includegraphics[width=\textwidth]{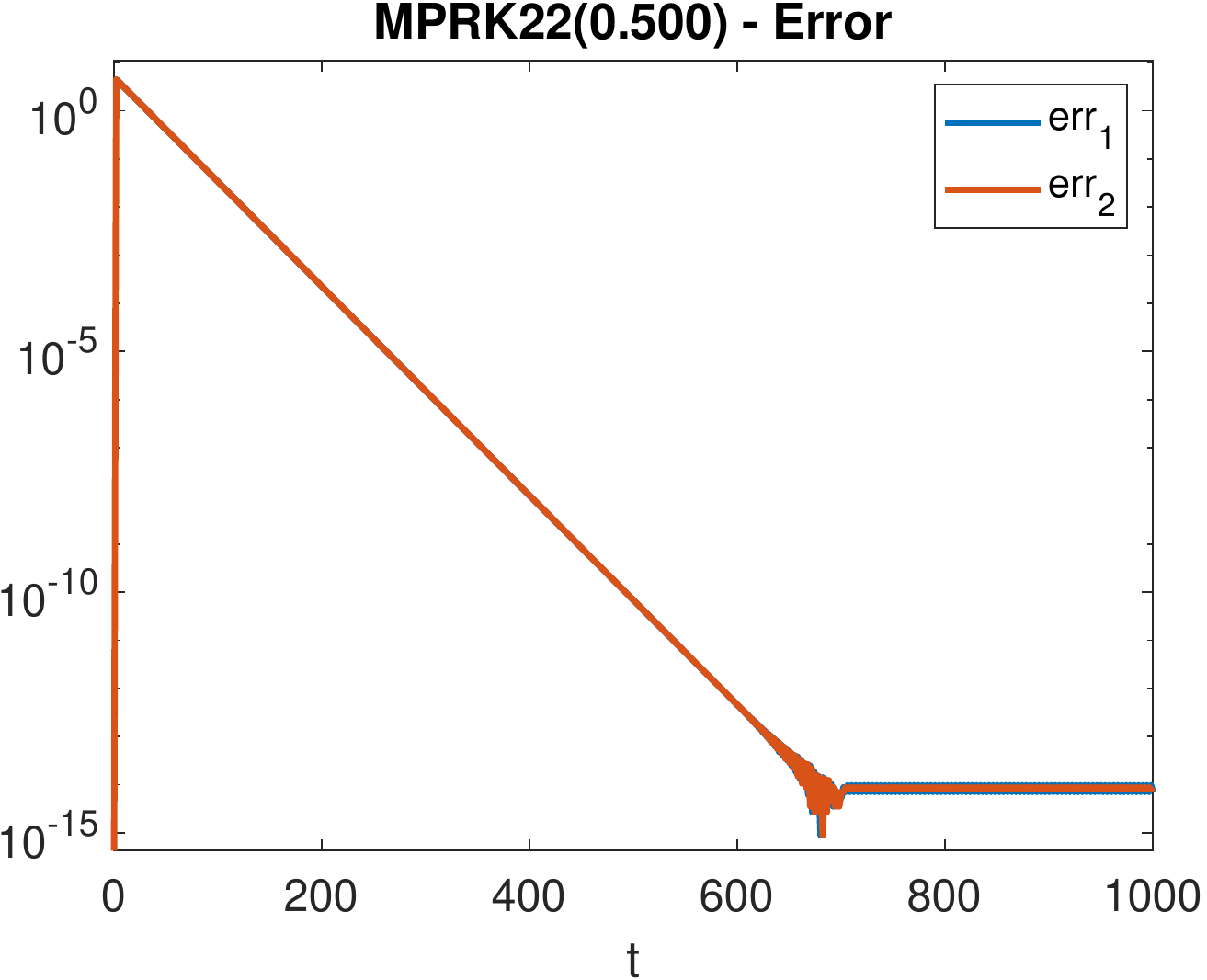}
		\subcaption{}\label{subfig:err0.5}
	\end{subfigure}\\
	\begin{subfigure}[t]{0.4\textwidth}
		\includegraphics[width=\textwidth]{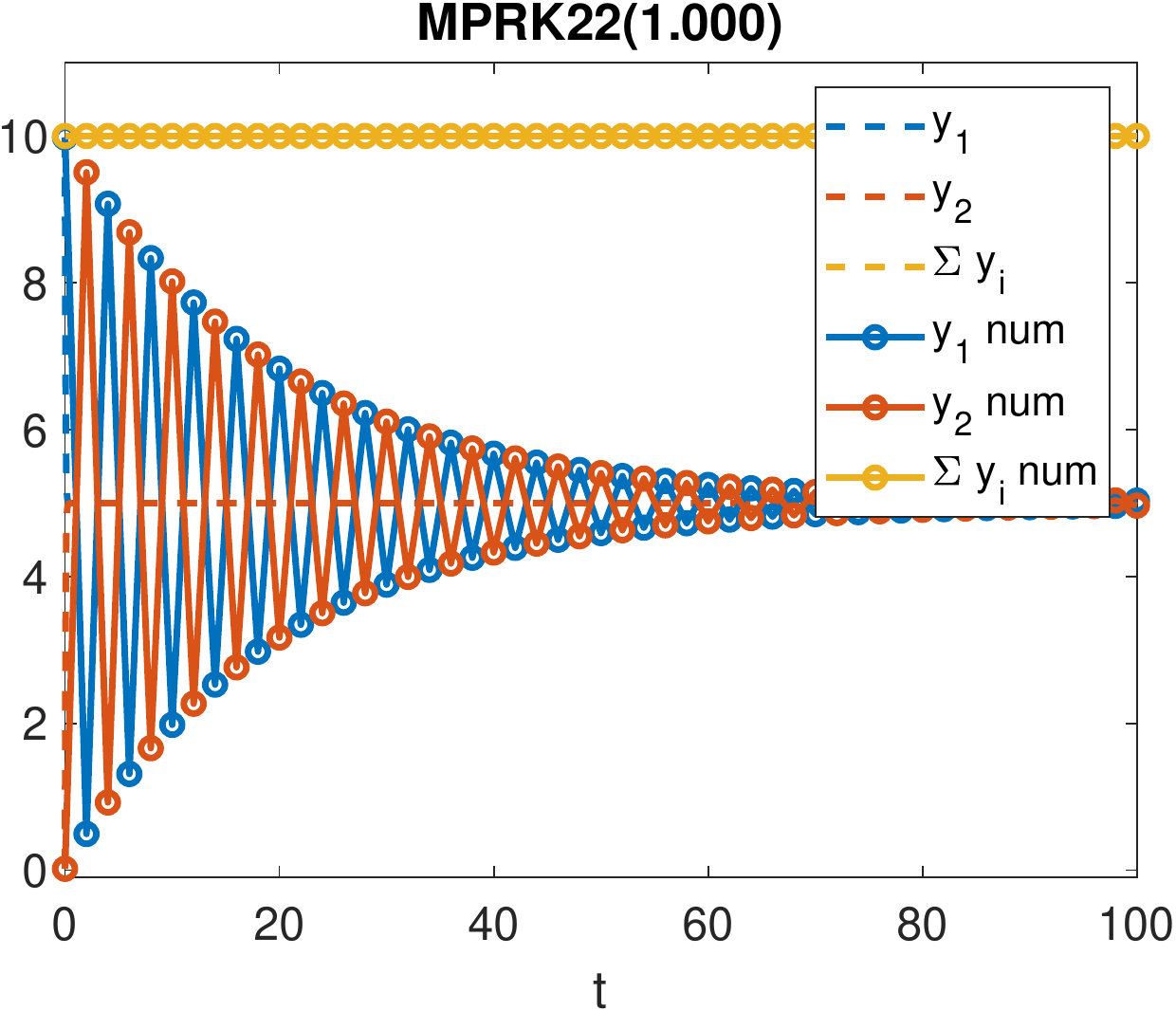}
		\subcaption{}\label{subfig:1}
	\end{subfigure}
	\begin{subfigure}[t]{0.423\textwidth}
		\includegraphics[width=\textwidth]{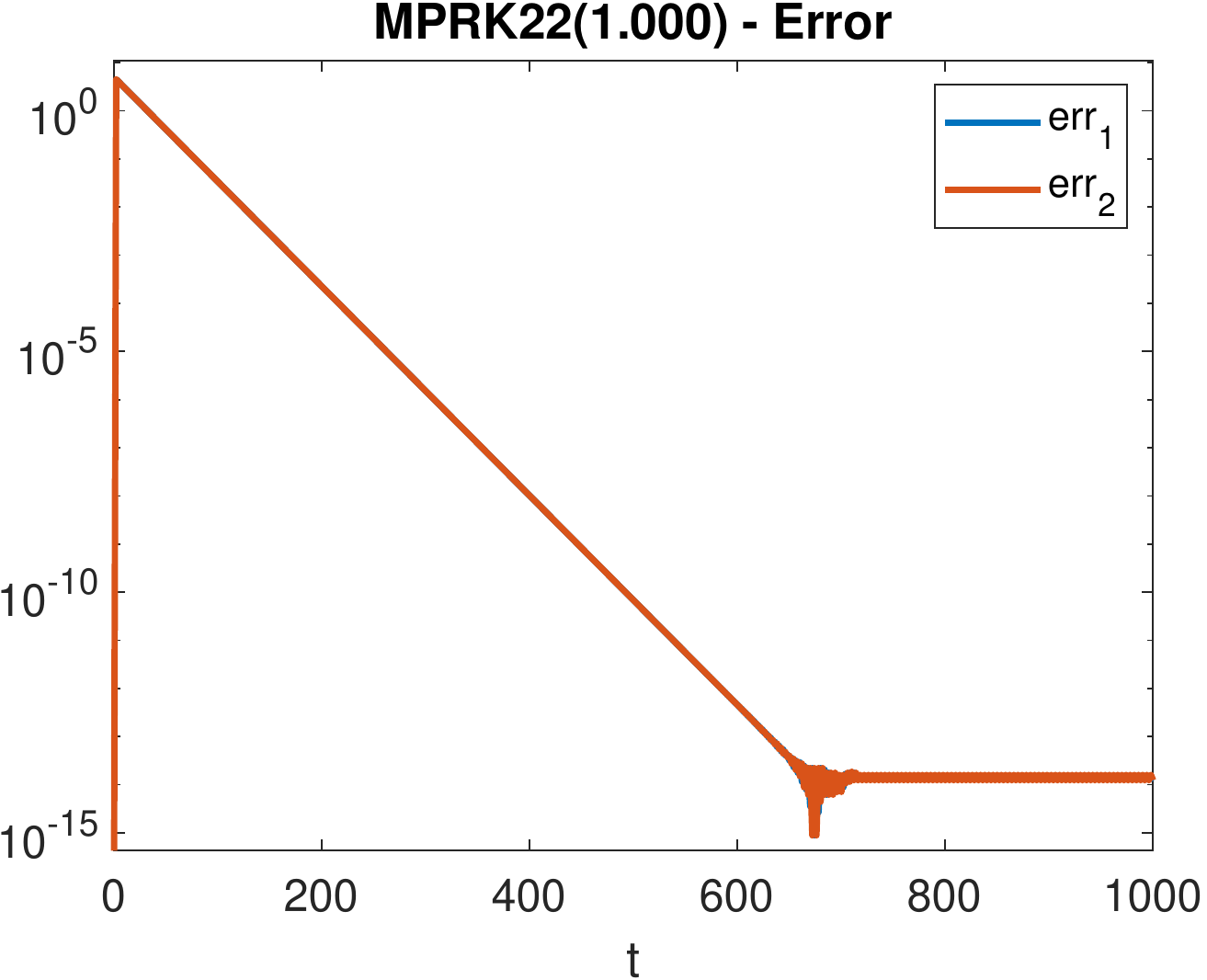}
		\subcaption{}\label{subfig:err1}
	\end{subfigure}\\
	\begin{subfigure}[t]{0.4\textwidth}
		\includegraphics[width=\textwidth]{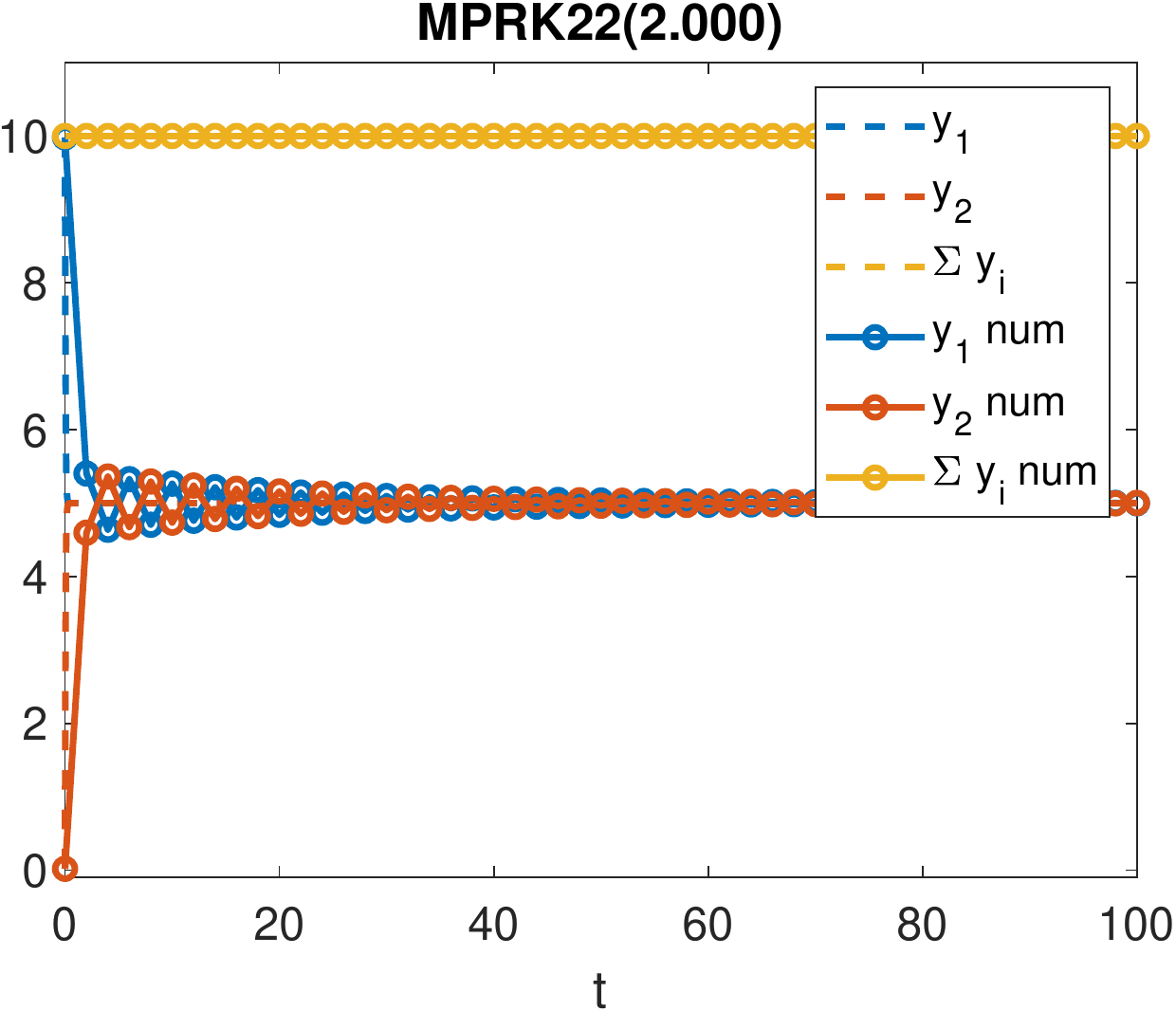}
		\subcaption{}\label{subfig:2}
	\end{subfigure}
	\begin{subfigure}[t]{0.423\textwidth}
		\includegraphics[width=\textwidth]{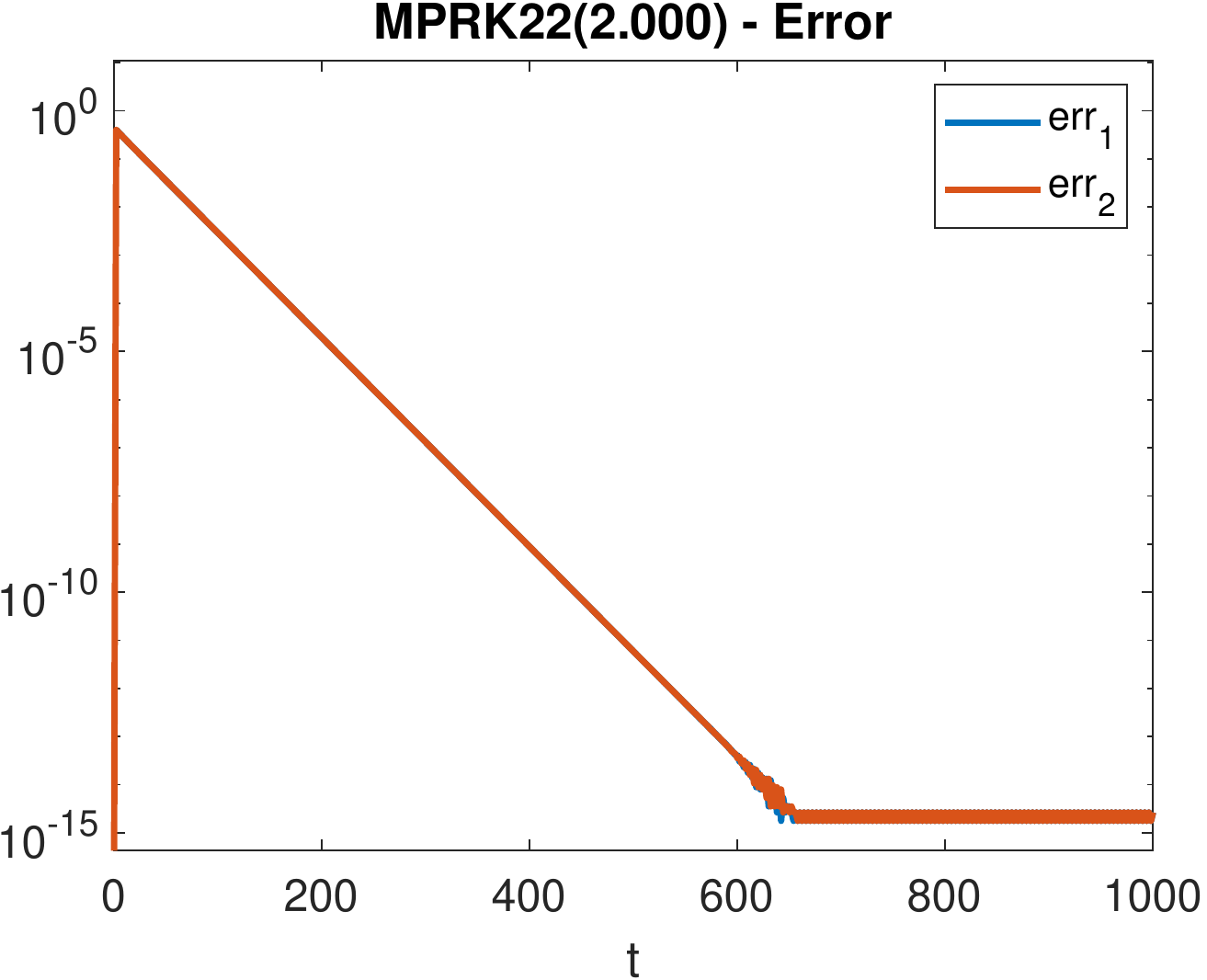}
		\subcaption{}\label{subfig:err2}
	\end{subfigure}
	\caption{Numerical approximations of \eqref{eq:PDS_test} using MPRK22($\alpha$) schemes for $\alpha\in\{0.5,1,2\}$, see Figures \ref{subfig:0.5}, \ref{subfig:1} and \ref{subfig:2}. There, the dashed lines indicate the exact solution \eqref{eq:anasol}. In the Figures \ref{subfig:err0.5}, \ref{subfig:err1} and \ref{subfig:err2} the corresponding errors $\operatorname{err}_i=\abs{y_i^n-y_i(n\Delta t)}$ are plotted for $i=1,2$.}\label{Fig:exp}
\end{figure}

In this section we numerically validate the statements of Corollary \ref{Cor:MPRKstab} and Corollary \ref{Cor:MPRKstab1}. For this, we solve the test problem \eqref{eq:PDS_test} with initial condition $\by^0=(9.98,0.02)^\intercal$. The analytical solution is plotted in Figure \ref{fig:analytic}, where one can see that the steady state solution is approximately reached at time $t=0.4$.
For the numerical approximation we use a comparably large time step size of $\Delta t=2$ to test the statements of Corollary \ref{Cor:MPRKstab} and Corollary \ref{Cor:MPRKstab1}. Indeed, Figure \ref{Fig:exp} reflects the theoretical results very well. 
Numerical approximations can be seen in Figures \ref{subfig:0.5}, \ref{subfig:1} and \ref{subfig:2}.
The corresponding error plots are depicted in Figures \ref{subfig:err0.5}, \ref{subfig:err1} and \ref{subfig:err2} showing that the error is of magnitude $10^{-14}$ for $t\geq 700$.

\section{Summary and Conclusion}
In this work we used Theorem \ref{Thm_MPRK_stabil} for the stability analysis of MPRK22($\alpha$) when applied to the nonlinear test problem \eqref{eq:PDS_test}. This is the first time that the stability theory of \cite{IKM22} is used in the context of a nonlinear differential equation. We discovered that all steady states of \eqref{eq:PDS_test} are stable fixed points of the method regardless of the chosen time step size. Furthermore, we proved that the iterates of MPRK22($\alpha$) locally converge towards the correct steady state solution of the initial value problem \eqref{eq:PDS_test}, \eqref{eq:IC}. The theoretical results are numerically validated and underline the robustness of the MPRK22($\alpha$) schemes.
\section{Acknowledgements}
The author Th.\ Izgin gratefully acknowledges the financial support by the Deutsche Forschungsgemeinschaft (DFG) through grant ME 1889/10-1.

		

		
	\end{document}